\documentclass[a4paper,12pt]{article}
\usepackage{amssymb}

\newcommand{\R}{\mathbb{R}}

\newcommand{\N}{\mathbb{N}}

\newcommand{\cuad}{{\sqcap\kern-.68em\sqcup}}

\newcommand{\norm}[1]{\|#1\|}

\newtheorem{definition}{Definition}[section]
\newtheorem{theorem}{Theorem}[section]
\newtheorem{proposition}{Proposition}[section]

\newtheorem{lemma}{Lemma}[section]

\newtheorem{remark}{Remark}[section]
\newcommand{\bremark}{\begin{remark} \em}
\newcommand{\eremark}{\end{remark} }

\begin{document}	
\begin{center}
{\Large \bf

Boundary blow-up solutions of elliptic  equations

\medskip

 involving regional fractional Laplacian

 }

  \bigskip  \medskip

 Huyuan Chen\footnote{chenhuyuan@yeah.net}  \qquad  Hichem Hajaiej\footnote{hh62@nyu.edu}

\smallskip

\begin{abstract}

In this paper, we study existence of boundary blow-up solutions for
elliptic equations involving  regional fractional Laplacian:
\begin{equation}\label{0}
  \arraycolsep=1pt
\begin{array}{lll}
 \displaystyle  (-\Delta)^\alpha_\Omega   u+ f(u)=0\qquad & {\rm in}\quad   \Omega,\\[2mm]
\phantom{ (-\Delta)^\alpha   +f(u) }
 \displaystyle   u=+\infty\quad & {\rm on}\quad   \partial  \Omega,
\end{array}
\end{equation}
where $\Omega$ is a bounded open domain in $\R^N$ ($N\ge 2$) with $C^2$ boundary $\partial\Omega$,
 $\alpha\in(0,1)$ and the operator $(-\Delta)^\alpha_\Omega$ is the regional fractional Laplacian.
When $f$ is a nondecreasing continuous function satisfying $f(0)\ge0$  and some additional conditions,
 we address the existence and nonexistence of solutions for problem  (\ref{0}). Moreover,  we further analyze the
 asymptotic behavior of solutions to problem  (\ref{0}).

\end{abstract}

\end{center}
  \vspace{1mm}


\vspace{2mm}

\setcounter{equation}{0}
\section{Introduction}

 The usual Laplaciain operator may be thought as a macroscopic manifestation of the Brownian motion, as known from
the Fokker-Plank equation for a stochastic
differential equation with a  Brownian motion (a Gaussian process), whereas the fractional Laplacian
operator $(-\Delta)^{ \alpha }$  is associated with
  a $2\alpha$-stable L\'{e}vy motion
(a non-Gaussian process) $L_t^{2\alpha}$, $\alpha\in(0,1)$, (see \cite{Duan2015} for a discussion about this microscopic-macroscopic relation.)
Given a  bounded open domain $\Omega$  in $\R^N$, the regional fractional Laplacians defined in $\Omega$ are generators of the reflected symmetric $2\alpha$-stable processes,  see \cite{CK,CKS,GM}.
Motivated by numerous applications related to (\ref{0}) and by the great mathematical interest in solving (\ref{0}) itself,
we tackle this rich PDE problem in this paper.

Let  $\Omega$ be a bounded open domain in $\R^N$ ($N\ge 2$) with $C^2$ boundary $\partial\Omega$, $\rho(x)={\rm dist}(x,\R^N\setminus\Omega)$
and  $f:\R\to\R$ be a nondecreasing, locally Lipschitz  continuous function satisfying $f(0)\ge0$.
We are concerned with the existence of boundary blow-up solutions for
elliptic equations involving  regional fractional Laplacian
\begin{equation}\label{eq 1.1}
  \arraycolsep=1pt
\begin{array}{lll}
 \displaystyle  (-\Delta)^\alpha_\Omega   u+ f(u)=0\qquad & {\rm in}\quad   \Omega,\\[2mm]
\phantom{ (-\Delta)^\alpha   +f(u)  }
 \displaystyle   u=+\infty\quad & {\rm on}\quad   \partial  \Omega,
\end{array}
\end{equation}
where   $\alpha\in(0,1)$ and  $(-\Delta)^\alpha_\Omega$ is the regional fractional Laplacian
 defined by
$$ (-\Delta)^\alpha_\Omega  u(x)=P.V. \int_{\Omega}\frac{u(x)-u(y)}{|x-y|^{N+2\alpha}}dy,\quad \ \ x\in\Omega.
$$
Here $P.V.$ denotes the principal value of the integral, that for notational simplicity we omit in what follows.

When $\alpha=1$, in  the seminal works by  Keller \cite{K} and Osserman \cite{O}, the authors  studied the boundary blow-up solutions for
 the nonlinear reaction diffusion equation
\begin{equation}\label{eq 12-1}
  \arraycolsep=1pt
\begin{array}{lll}
 -\Delta u+f(u)=0\quad &
 \rm{in}\ \ \ &\Omega,\\[2mm]
 \phantom{----- }
u=+\infty\quad&  \rm{on}\ \ \ &\partial\Omega.
\end{array}
\end{equation}
They independently proved that this equation admits a solution if and only if  $f$ is a nondecreasing positive function satisfying the Keller-Osserman criterion, that is,
\begin{equation}\label{ko}
\int_1^{+\infty}\frac{ds}{\sqrt{\int_0^{s} f(t)dt}}<+\infty.
\end{equation}
From then on, boundary blow-up problem (\ref{eq 12-1}) has been extended by numerous mathematicians in various ways: weakening
the assumptions on the domain, generalizing the differential operator and the nonlinear term for equations and systems.
Moreover, the qualitative properties of   boundary blow-up solutions, such as asymptotic behavior, uniqueness and symmetry results,
attract a great attention, see the references \cite{AGQ,BM1,BM,FQ0,MV0,MV}.

In a recent work,  Chen-Felmer-Quaas \cite{CFQ}  considered an analog of (\ref{eq 12-1}) where the Laplacian is replaced by the fractional Laplacian
 \begin{equation}\label{eq 12-3}
 \arraycolsep=1pt
\begin{array}{lll}
  (-\Delta)^\alpha u+f(u)=0\quad&
 {\rm in}\quad \Omega,\\[1.5mm]
  \phantom{------\ }
  u=0 & {\rm in}\quad \R^N\setminus\Omega,\\[1.5mm]
  \phantom{\ \ }
\displaystyle \lim_{x\in\Omega,x\to\partial\Omega}u(x)=+\infty,
\end{array}
\end{equation}
where the fractional Laplacian operator $(-\Delta)^\alpha$ is defined as
$$ (-\Delta)^\alpha  u(x)=P.V. \int_{\R^N}\frac{u(x)-u(y)}{|x-y|^{N+2\alpha}}dy.
$$
They studied the existence, uniqueness and non-existence of boundary blow-up solutions by Perron's method when
$f(s)=s^p$ with $p>1$. Later on, the authors and Wang in  \cite{CHY1} studied the boundary blow-up solutions of (\ref{eq 12-3}) which is derived by measure type data when $f$ is a continuous
and  increasing function satisfying
\begin{equation}\label{1.1}
 \int_1^{\infty} f(s)s^{-1-\frac{1+\alpha}{1-\alpha}}ds<+\infty.
\end{equation}
We obtained a sequence of boundary blow-up solutions of (\ref{eq 12-3}), which  have the asymptotic  behavior
${\rm dist}(x,\partial\Omega)^{\alpha-1}$ as  $x\to\partial\Omega$. In particular, when $f(s)\le c_1s^q$  for $s\ge0$, where  $q\le 2\alpha+1$ and $c_1>0$, this sequence of
solutions blow up every where in $\Omega$.

For a regular function $u$ such that $u=0$ in $\R^N\setminus\bar\Omega$, we remark that
$$
   (-\Delta)^\alpha_\Omega u(x)= (-\Delta)^\alpha u(x)-u(x)\phi(x),\quad\forall x\in\Omega,
$$
where
$$
 \phi(x)=\int_{\R^N\setminus \Omega} \frac1{|x-y|^{N+2\alpha}}dy.
$$
From the connections between the fractional Laplacian and the regional fractional Laplacian, we observe that the boundary blowing up solution of
 (\ref{eq 12-3}) provides a sub solution for (\ref{eq 1.1}),  then  we have following proposition.

\begin{proposition}\label{cr 1.1}
Assume that $\alpha\in(0,1)$  and    $f$ is a nondecreasing function satisfying $f(0)\ge0$ and locally Lipschitz continuous in $\R$ .

$(i)$ If $f(s)\le c_1s^q$  for $s\ge0$, where  $q\le 2\alpha+1$ and $c_1>0$, then problem (\ref{eq 1.1}) has no  solution $u$ satisfying
\begin{equation}\label{1.3}
 \lim_{\rho(x)\to0^+}u(x)\rho(x)^{1-\alpha}=+\infty.
\end{equation}

$(ii)$ If
\begin{equation}\label{1.5}
 c_2s^p\le  f(s)\le c_3s^q\quad{\rm  for}\quad s\ge1,
\end{equation}
 where  $2\alpha+1< p\le q \le \frac{1+\alpha}{1-\alpha}$ and $c_2,c_3>0$, then problem (\ref{eq 1.1}) has a solution $u$ satisfying
\begin{equation}\label{1.4}
c_4 \rho(x)^{-\frac{2\alpha}{q-1}}\le  u(x)\le c_5 \rho(x)^{-\frac{2\alpha}{p-1}},\quad \forall x\in\Omega,
\end{equation}
where $c_5\geq c_4>0$.
\end{proposition}

We notice that Proposition \ref{cr 1.1} can not cover the case where $f(s)\ge s^p$ with $p\ge \frac{1+\alpha}{1-\alpha}$. Our purpose in this note
is to solve more general cases. To this end,  we first introduce an important proposition on the  regional fractional elliptic problem with finite boundary data.

\begin{proposition}\label{pr 3.1}
Let  $\alpha\in(\frac12,1)$, $n\in\N$,  $g\in C^1(\bar \Omega)$ and $f$ be a locally Lipschitz continuous and nondecreasing function.

Then problem
\begin{equation}\label{eq 3.1}
  \arraycolsep=1pt
\begin{array}{lll}
 \displaystyle  (-\Delta)^\alpha_\Omega   u+f(u)=g \qquad & {\rm in}\quad   \Omega,\\[2mm]
\phantom{ (-\Delta)^\alpha +f(u) }
 \displaystyle   u=n\quad & {\rm on}\quad   \partial  \Omega
\end{array}
\end{equation}
 admits a unique solution $u_n$ such that
\begin{equation}\label{3.1}
-c_6\left(\norm{g_-}_{L^\infty(\Omega)}+f(n)\right)\rho^{2\alpha-1}\le  u_n-n\le c_6\norm{g_+}_{L^\infty(\Omega)}\rho^{2\alpha-1}\quad{\rm in}\quad \Omega,
\end{equation}
 where  $g_\pm=\max\{\pm g,0\}$ and  $c_6>0$ is independent of $n$, $f$ and $g$.

Moreover, if $g\ge 0$  and $f(0)\ge 0$, then $u_n$ is positive.
\end{proposition}

The derivation of the solution of (\ref{eq 3.1}) makes use of the Green's function of the regional fractional Laplacian  and  Perron's method. The authors in \cite{CK} showed that for $\alpha\in (\frac12,1)$, the  Green's function of the regional fractional provides boundary decay estimate, while for $\alpha\in (0,\frac12]$,
the  Green's function of the regional fractional behaviors very different, without any boundary decaying, thus it is even hard to obtain a solution for (\ref{eq 3.1}).

We call a solution $u_m$ of (\ref{eq 1.1}) is the minimal solution if for any solution $v$ of (\ref{eq 1.1}), we have that
$$v\ge u_m\quad {\rm in}\quad \Omega.$$
As normal, the minimal boundary blow-up solution of with $\alpha\in (\frac12,1)$ is approached by the solutions of (\ref{eq 3.1}) by  taking $n\to+\infty$.

\begin{theorem}\label{teo 1}
Assume that  $\alpha\in(\frac12,1)$  and  $f$ is a nondecreasing continuous function satisfying $f(0)\ge0$. Furthermore,

$(i)$ If $f(s)\ge c_7s^p$ for $s\ge0$, where  $p>1+2\alpha$ and $c_7>0$, then  problem (\ref{eq 1.1})  possesses the minimal boundary blow-up solution $u_m$.

Assume more that
$f(s)\le c_8s^q$ for $s\ge1$, where $q\ge p$ and $c_8>0$, then $u_m$
has asymptotic behavior near the boundary as
\begin{equation}\label{1.2}
 c_9\rho(x)^{-\frac{2\alpha-1}{q-1}}\le u_m(x)\le c_{10}\rho(x)^{-\frac{2\alpha}{p-1}},
\end{equation}
where $c_{10}\geq c_9>0$.

$(ii)$   If  $f(s)\le c_{11}s^q$  for $s\ge0$, where $c_{11}>0$ and
\begin{equation}\label{q+1}
 q\le 1+2\alpha \quad{\rm and}\quad q< \frac{\alpha}{1-\alpha},
\end{equation}
 then problem (\ref{eq 1.1}) has no solution.

\end{theorem}

Compared to Proposition \ref{cr 1.1}, we notice that $(i)$ when $\alpha\in(\frac12,1)$, we improve the existence for the case that
  $f(s)\ge c_7s^p$ for $s\ge0$ and  $p>1+2\alpha$ in  Theorem \ref{teo 1};
  $(ii)$ if $\alpha>\frac{\sqrt{2}}{2}$ for $f(s)=s^p$ with $p\le 1+2\alpha$, problem (\ref{eq 1.1}) has  any solution.

  The lower bound in (\ref{1.2}) is derived by the  inequality  (\ref{3.1}) and the upper bound in (\ref{1.2}) is obtained by constructing a suitable super-solution
for problem (\ref{eq 1.1}).

 This article is organized as follows.  Section \S2 is devoted to present some preliminaries on the definition of viscosity solution, Comparison Principle, Stability theorem,  regularity results and to make use of solutions of corresponding problem with the fractional Laplacian to prove Proposition \ref{cr 1.1}.   In Section \S3,  we first prove the existence of solutions in order to problem (\ref{3.1}), asymptotic behavior
   and then  prove Theorem \ref{teo 1}.

\setcounter{equation}{0}
\section{Preliminary}

The purpose of this section is to introduce some preliminaries. We start it by defining the notion of viscosity solution, inspired by the definition of
viscosity sense for nonlocal problem in \cite{CS}.

\begin{definition}\label{de 2.2}
We say that a  continuous function $u\in L^1(\Omega)$ is a viscosity super-solution (sub-solution)
 of
 \begin{equation}\label{eq 2.1}
  \arraycolsep=1pt
\begin{array}{lll}
 \displaystyle  (-\Delta)^{\alpha}_\Omega u+f(u)=g\qquad & {\rm in}\quad   \Omega,\\[2mm]
\phantom{ (-\Delta)^\alpha   +f(u)}
 \displaystyle   u=h\quad & {\rm on}\quad   \partial  \Omega,
\end{array}
\end{equation}
if $u\ge h$ (resp. $u\le h$) on $\partial  \Omega$ and
for every point $x_0\in\Omega$ and some  neighborhood $V$ of
$x_0$ with $\bar V\subset \Omega $ and for any $\varphi \in
C^2(\bar V)$ such that $u(x_0)=\varphi(x_0)$ and $x_0$ is the minimum (resp. maximum) point of $u-\varphi$ in $V$,  let
\begin{eqnarray*}
\tilde u  =\left\{ \arraycolsep=1pt
\begin{array}{lll}
\varphi\ \ \ & \rm{in}\ \ &V,\\[2mm]
u \ \ & \rm{in}\ \ &\Omega\setminus V,
\end{array}
\right.
\end{eqnarray*}
we have
$$(-\Delta)^{\alpha }_\Omega\tilde u(x_0)+f(  u(x_0))\ge g(x_0)\quad (resp.  \ (-\Delta)^{\alpha }_\Omega\tilde u(x_0)+f(  u(x_0))\le g(x_0)).$$

We say that $u$ is a
viscosity solution of (\ref{eq 2.1}) if it is  a viscosity super-solution and also a viscosity sub-solution of (\ref{eq 2.1}).
\end{definition}

Now we introduce the Comparison Principle.

\begin{theorem}\label{comparison}
Assume that the functions  $g:\Omega\to\R$, $h:\partial\Omega\to\R$ are
continuous and $f:\R\to\R$ is nondecreasing. Let $u$ and $v$ be a viscosity super-solution and
sub-solution  of (\ref{eq 2.1}), respectively.
If $$v \le u \quad{\rm on}\quad    \partial\Omega,$$
 then
\begin{equation}\label{2.1}
 v  \le u \quad{\rm in}\quad  \Omega.
\end{equation}

\end{theorem}
{\bf Proof.}  Let us define $w=u-v$, then
\begin{equation}\label{eq 2.2}
  \arraycolsep=1pt
\begin{array}{lll}
 \displaystyle  (-\Delta)^{\alpha}_\Omega w \ge f(v)-f(u)\qquad & {\rm in}\quad   \Omega,\\[2mm]
\phantom{ (-\Delta)^\alpha   }
 \displaystyle  w\ge 0\quad & {\rm on}\quad   \partial  \Omega.
\end{array}
\end{equation}
If (\ref{2.1}) fails, then there exists $x_0\in\Omega$ such that
$$w(x_0)=u(x_0)-v(x_0)=\min_{x\in\Omega}w(x)<0,$$
by the fact that $f$ is nondecreasing, we have that $f(v(x_0))-f(u(x_0))\ge 0$  and then in the viscosity sense,
\begin{equation}\label{y 2.1}
 (-\Delta)^{\alpha}_\Omega w(x_0)\ge 0.
\end{equation}
 Since $w$ is a viscosity super solution  $x_0$ is the minimum point in $\Omega$ and $w\ge 0$ on $\partial\Omega$, then we can
take a small neighborhood $V_0$ of $x_0$ such that $\tilde w=w(x_0)$ in $V_0$,
From (\ref{y 2.1}), we have that
$$ (-\Delta)^{\alpha}_\Omega \tilde  w(x_0)\ge 0.$$
But
$$(-\Delta)^{\alpha}_\Omega \tilde w(x_0)=   \int_{\Omega\setminus V_0} \frac{w(x_0)-w(y)}{|x_0-y|^{N+2\alpha}}dy<0,$$
which is impossible.  \qquad $\Box$

\smallskip
For a regular function $w$ such that $w=0$ in $\R^N\setminus\bar\Omega$, we observe that
  \begin{equation}\label{con}
   (-\Delta)^\alpha_\Omega w(x)= (-\Delta)^\alpha w(x)-w(x)\phi(x),\quad\forall x\in\Omega,
  \end{equation}
where
\begin{equation}\label{con 1}
 \phi(x)=\int_{\R^N\setminus \Omega} \frac1{|x-y|^{N+2\alpha}}dy.
\end{equation}
 \begin{lemma}\label{lm 2.2}
Let $\phi$ be defined in (\ref{con 1}) and $\rho(x)=dist(x,\partial\Omega)$, then  $\phi\in C^{0,1}_{\rm loc}(\Omega)$ and
\begin{equation}\label{2.6}
 \frac1{c_{12}}\rho(x)^{-2\alpha}\le \phi(x)\le c_{12}\rho(x)^{-2\alpha},\quad x\in\Omega,
\end{equation}
for some $c_{12}>0$.
\end{lemma}
{\bf Proof.} For
 $x_1,x_2\in \Omega $ and  any $z\in \R^N\setminus \Omega$, we have that
 $$|z-x_1|\ge \rho(x_1)+\rho(z), \qquad |z-x_2|\ge \rho(x_2)+\rho(z)$$ and
$$||z-x_1|^{N+2\alpha}-|z-x_2|^{N+2\alpha}|\leq
c_{13}|x_1-x_2|(|z-x_1|^{N+2\alpha-1}+|z-x_2|^{N+2\alpha-1}),
$$
for some $c_9>0$ independent of $x_1$ and $x_2$. Then
\begin{eqnarray*}
&&|\phi(x_1)-\phi(x_2)| \leq\int_{\R^N\setminus \Omega}  \frac{||z-x_2|^{N+2\alpha}-|z-x_1|^{N+2\alpha}|}{|z-x_1|^{N+2\alpha}|z-x_2|^{N+2\alpha}}dz
\\&&\leq c_{13}|x_1-x_2|\left[\int_{\R^N\setminus \Omega}  \frac{dz}{|z-x_1||z-x_2|^{N+2\alpha}}+\int_{\R^N\setminus \Omega}  \frac{dz}{|z-x_1|^{N+2\alpha}|z-x_2|}\right].
\end{eqnarray*}
By direct computation, we have that
\begin{eqnarray*}
 \int_{\R^N\setminus \Omega}  \frac{1}{|z-x_1||z-x_2|^{N+2\alpha}}dz
  &\le& \int_{\R^N\setminus{B_{\rho(x_1)}(x_1)}} \frac{1}{|z-x_1|^{N+2\alpha+1}}dz
 \\&&+\int_{\R^N\setminus{B_{\rho(x_2)}(x_2)}} \frac{1}{|z-x_2|^{N+2\alpha+1}}dz \\
 &\le& c_{14}[\rho(x_1)^{-1-2\alpha}+\rho(x_2)^{-1-2\alpha}]
\end{eqnarray*}
and similar to obtain that
$$\int_{\R^N\setminus \Omega} \frac{1}{|z-x_1|^{N+2\alpha}|z-x_2|}dz\le  c_{14}[\rho(x_1)^{-1-2\alpha}+\rho(x_2)^{-1-2\alpha}],$$
where $c_{14}>0$ is independent of $x_1,x_2$.
Then
$$|\phi(x_1)-\phi(x_2)|\le c_{13} c_{14}[\rho(x_1)^{-1-2\alpha}+\rho(x_2)^{-1-2\alpha}]|x_1-x_2|,$$
that is,
 $\phi$ is $C^{0,1}$ locally in $\Omega$.

  Now we prove (\ref{2.6}). Without loss of generality, we may assume that
$0\in\partial\Omega$, the inside pointing normal vector at $0$ is $e_N=(0,\cdots,0,1)\in\R^N$ and 
let $s\in(0,\frac14)$ such that $\R^N\setminus \Omega\subset \R^N\setminus B_s(se_N)$ and for $c>0$, we denote the cone
$$A_c=\{y=(y',y_N)\in\R^N: y_N\le s-c|y'|\}.$$
We observe that
there is $c_{15}>0$ such that
$$\left[A_{c_{15}}\cap \left(B_1(se_N)\setminus B_{2s}(se_N)\right)\right]\subset \R^N\setminus \Omega. $$
By the definition of $\phi$, we have that
 \begin{eqnarray*}
\phi(se_N)= \int_{\R^N\setminus \Omega} \frac1{|se_N-y|^{N+2\alpha}}dy \le  \int_{\R^N\setminus B_s(se_N)} \frac1{|se_N-y|^{N+2\alpha}}dy  \le c_{16}s^{-2\alpha}
 \end{eqnarray*}
 for some $c_{16}>0$. On the other hand, we have that
\begin{eqnarray*}
 \int_{\R^N\setminus \Omega} \frac1{|se_N-y|^{N+2\alpha}}dy \ge  \int_{A_{c_{15}}\cap \left(B_1(se_N)\setminus B_{2s}(se_N)\right)} \frac1{|se_N-y|^{N+2\alpha}}dy
     \ge c_{17}s^{-2\alpha},
 \end{eqnarray*}
for some $c_{17}\in(0,1)$.
The proof ends.\qquad $\Box$

\smallskip

 The next theorem gives the stability property for viscosity solutions in our setting.

\begin{theorem}\label{stability}
Assume that the function  $g:\Omega\to\R$ is
continuous,  $f:\R\to\R$ is nondecreasing and $f(0)\ge0$.  Let $(u_n)_n$, $n\in\N$ be a sequence of
 functions in  $C^1(\Omega)$, uniformly  bounded in
$L^1(\Omega)$, $g_n$ and $g$
be continuous in $\Omega$ such that

$(-\Delta)^\alpha_\Omega u_n+ f(u_n)\ge g_n\ ({\rm resp.}\ (-\Delta)^\alpha_\Omega
u_n+ f(u_n)\le g_n)$ in $\Omega$ in viscosity sense,

$u_n\to u$  locally uniformly in  $\Omega$,

$ u_n\to u$  in  $L^1(\Omega)$,

$ h_n\to h$ locally uniformly in  $\Omega$.\\
    Then $(-\Delta)^\alpha_\Omega  u+f(u)\ge g\ ({\rm resp.}\ (-\Delta)^\alpha_\Omega u+f(u)\le g)$ in $\Omega$ in the viscosity sense.
\end{theorem}
{\bf Proof.}  We define $\tilde u_n=u_n$ in $\Omega$, $\tilde u_n=0$ in $\R^N\setminus\bar \Omega$ and
$\tilde u=u$ in $\Omega$, $\tilde u=0$ in $\R^N\setminus\bar \Omega$, then
$$(-\Delta)^\alpha_\Omega u_n(x)=(-\Delta)^\alpha \tilde u_n(x)-u_n(x)\phi(x),\qquad x\in\Omega.$$
where $\phi$ is defined as (\ref{con 1}).
By Lemma \ref{lm 2.2},
  $\phi\in C^{0,1}_{\rm loc}(\Omega)$ and
$ \phi(x)\le c_8\rho(x)^{-2\alpha},$ $x\in\Omega.$ They we apply
 \cite[Theorem 2.4]{CFQ} to obtain that
 $(-\Delta)^\alpha   \tilde u+f(\tilde u)\ge g+\phi \tilde u$ (resp. $(-\Delta)^\alpha  \tilde u+f(\tilde u)\le g+\phi \tilde u $) in $\Omega$ in viscosity sense,
which implies$(-\Delta)^\alpha_\Omega  u+f(u)\ge g\ ({\rm resp.}\ (-\Delta)^\alpha_\Omega u+f(u)\le g)$ in $\Omega$ in viscosity sense.  \qquad$\Box$

\smallskip

Next we have an interior regularity result. For simplicity, we denote by $C^t$ the space $C^{t_0,t-t_0}$ for $t\in(t_0,t_0+1)$, $t_0$ is a positive integer.
\begin{proposition}\label{pr 2.1}
Assume that $\alpha\in(\frac12,1)$, $g\in C_{\rm loc}^{\theta}(\Omega)$ with $\theta>0$, $w\in C_{\rm loc}^{2\alpha+\epsilon}(\Omega)\cap L^1(\Omega)$ with $\epsilon>0$ 
and $2\alpha+\epsilon$  not being an integer is a solution of
\begin{equation}\label{homo}
(-\Delta)^\alpha_\Omega w=g\quad {\rm in}\ \ \Omega.
\end{equation}
Let $\mathcal{O}_1, \mathcal{O}_2$ be  open $C^2$ sets such that
$$\bar \mathcal{O}_1\subset\mathcal{O}_2\subset  \bar \mathcal{O}_2\subset \Omega.$$
Then\\
(i) for any $\gamma\in(0,2\alpha)$ not an integer, there exists $c_{18}>0$ such that
\begin{equation}\label{2.2}
\norm{w}_{C^\gamma(\mathcal{O}_1)}\le c_{18}\left[\norm{w}_{L^\infty(\mathcal{O}_2)}+\norm{w}_{L^1(\Omega)}+\norm{g}_{L^\infty(\mathcal{O}_2)} \right];
\end{equation}
(ii) for any $\epsilon'\in(0,\min\{\theta,\epsilon\})$, $2\alpha+\epsilon'$ not an integer, there exists $c_{19}>0$ such that
\begin{equation}\label{2.3}
\norm{w}_{C^{2\alpha+\epsilon'}(\mathcal{O}_1)}\le c_{19}\left[\norm{w}_{L^\infty(\mathcal{O}_2)}+\norm{w}_{L^1(\Omega)}+\norm{g}_{c^\theta(\mathcal{O}_2) }\right].
\end{equation}

\end{proposition}
{\bf Proof.} Let $\tilde w=w$ in $\Omega$, $\tilde w=0$ in $\R^N\setminus\bar \Omega$,  we have that
\begin{eqnarray*}
 (-\Delta)^\alpha \tilde w(x)  =  (-\Delta)^\alpha_\Omega w(x)+ w(x)\phi(x),\quad \forall x\in\Omega,
\end{eqnarray*}
where $\phi$ is defined as (\ref{con 1}).
It follows by Lemma \ref{lm 2.2}, $\phi\in C^{0,1}_{\rm loc}(\Omega)$.
Combining with (\ref{homo}), we have that
\begin{eqnarray*}
 (-\Delta)^\alpha \tilde w(x)  = g(x)+  w(x)\phi(x),\quad \forall x\in\Omega.
\end{eqnarray*}
By \cite[Lemma 3.1]{CV3}, for any $\gamma\in(0,2\alpha)$, we have that
\begin{eqnarray*}
\norm{w}_{C^\gamma(\mathcal{O}_1)}&\le& c_{20}\left[\norm{w}_{L^\infty(\mathcal{O}_2)}+\norm{w}_{L^1(\Omega)}+\norm{g+w\phi}_{L^\infty(\mathcal{O}_2)} \right] \\
   &\le& c_{21}\left[\norm{w}_{L^\infty(\mathcal{O}_2)}+\norm{w}_{L^1(\Omega)}+\norm{g}_{L^\infty(\mathcal{O}_2)} \right]
\end{eqnarray*}
 and by  \cite[Lemma 2.10]{RS}, for any $\epsilon'\in(0,\min\{\theta,\epsilon\})$, we have that
 \begin{eqnarray*}
\norm{w}_{C^{2\alpha+\epsilon'}(\mathcal{O}_1)}&\le& c_{22}\left[\norm{w}_{C^{\epsilon'}(\mathcal{O}_2)}+\norm{g+w\phi}_{C^{\epsilon'}(\mathcal{O}_2)} \right] \\
   &\le& c_{23}\left[\norm{w}_{L^\infty(\mathcal{O}_2)}+\norm{w}_{L^1(\Omega)}+\norm{g}_{C^{\epsilon'}(\mathcal{O}_2)} \right],
\end{eqnarray*}
where $c_{22},c_{23}>0$. This ends the proof.\qquad$\Box$

\smallskip

\subsection{Proof of Proposition \ref{cr 1.1}}

Basically, the existence for boundary blow-up problem is usually   resorted to the Perron's  method. In this subsection,
we extend the Perron's method to the problem involving regional fractional Laplacian.

To this end, we first introduce the existence of boundary blow-up solution of fractional elliptic problem with locally Lipschitz continuous nonlinearity $f$,
precisely,
\begin{equation}\label{a.1.1}
\left\{ \arraycolsep=1pt
\begin{array}{lll}
 (-\Delta)^{\alpha} u(x)+f(u)=g,\ \ \ \ &
x\in\Omega,\\[2mm]
u(x)=0,& x\in\bar\Omega^c,\\[2mm]
\lim_{x\in\Omega,\ x\to\partial\Omega}u(x)=+\infty.
\end{array}
\right.
\end{equation}

\begin{theorem}\label{th 3.2.4-27}
  Assume that $f:\R\to \R$ is nondecreasing, $C^\gamma_{loc}$  and $f(0)=0$, the function $g:\Omega \to \R$ is a $C^\gamma_{loc}$ in $\Omega$.   Suppose that there
are super-solution $\bar U$ and sub-solution $\underline U$ of
(\ref{a.1.1}) such that $\bar U$ and $\underline U$ are $C^2$
locally in $\Omega$, bounded in
$L^1(\R^N,\frac{dy}{1+|y|^{N+2\alpha}})$ and
$$\bar U\geq \underline U \ \ {\rm in}\ \Omega,\ \ \ \
\liminf_{x\in\Omega,x\to\partial\Omega}\underline U(x)=+\infty, \ \ \ \
 \bar U= \underline U=0\ \ {\rm in}\ \bar\Omega^c.$$ Then
there exists at least one solution $u$ of (\ref{a.1.1}) in the
viscosity sense and
$$\underline U\leq u\leq\bar U \ \ \rm{in}\ \ \Omega.$$

Additionally, suppose that $g\ge0$ in $\Omega,$ then
$u>0$ in $\Omega$.
\end{theorem}
{\bf Proof.} We follow the  proof of \cite[Theorem 2.6]{CFQ} replacing $|u|^{p-1}u$ by $f(u)$.

\begin{theorem}\label{th 3.2.4}
Let $\Omega$ be an open  bounded $C^2$  domain and $p>0$.   Suppose that there
are super-solution $\bar U$ and sub-solution $\underline U$ of
(\ref{eq 1.1}) such that $\bar U$ and $\underline U$ are $C^2$
locally in $\Omega$,
$$\bar U\geq \underline U\ \ {\rm in}\ \Omega, \ \ \ \
\liminf_{x\in\Omega,x\to\partial\Omega}\underline U(x)=+\infty.$$ Then
there exists at least one solution $u$ of (\ref{eq 1.1}) in the
viscosity sense and
\begin{equation}\label{2.11}
 \underline U\leq u\leq\bar U \ \ \rm{in}\ \ \Omega.
\end{equation}

\end{theorem}
{\bf Proof.}  From (\ref{con}), to search the solution of (\ref{eq 1.1}) is equivalent to find out the solution of the fractional problem
\begin{equation}\label{eq 2.4}
 \arraycolsep=1pt
\begin{array}{lll}
  (-\Delta)^\alpha u+f(u)=\phi u\quad&
 {\rm in}\quad \Omega,\\[1.5mm]
  \phantom{------\ }
  u=0 & {\rm in}\quad \R^N\setminus\Omega,\\[1.5mm]
  \phantom{\ \ }
\displaystyle \lim_{x\in\Omega,x\to\partial\Omega}u(x)=+\infty,
\end{array}
\end{equation}
where $\phi$ is given by (\ref{con 1}). Make zero extensions of $\bar U$ and $\underline U$ in $\R^N\setminus\Omega$ and still denote them by $\bar U$ and $\underline U$
respectively, then $\bar U$ and $\underline U$ are the super and sub solutions of (\ref{eq 2.4}). Now we apply Theorem \ref{th 3.2.4-27}  to obtain the existence
of solution to (\ref{eq 2.4})

 From Lemma \ref{lm 2.2}, $\phi$ is $C^{0,1}$ locally in $\Omega$, so is $\phi \underline U$, then by Theorem \ref{th 3.2.4-27},
 we obtain that problem (\ref{eq 2.4}) replaced $\phi u$
 by $\phi \underline U$ admits a solution $u_1$ satisfying (\ref{2.11}). By regularity results in \cite{RS}, we have that
 $$\norm{u_1}_{C^{2\alpha+\gamma}(\Omega)}\le c_{24}\norm{\bar U}_{L^\infty(\Omega)}   $$
 for some $\gamma\in(0,1)$.

 Inductively, by Theorem \ref{th 3.2.4-27}, we obtain that problem (\ref{eq 2.4}) replaced $\phi u$
 by $\phi u_{n-1}$ has a solution $u_n$ such that
 \begin{equation}\label{2.12}
u_{n-1}\leq u_n\leq\bar U \ \ \rm{in}\ \ \Omega.
\end{equation}
We apply stability Theorem \cite[Theorem 2.4]{CFQ} and regularity result in \cite{RS}, we obtain that the limit of $\{u_n\}_n$ is a solution of (\ref{eq 2.4}).
 \qquad$\Box$

\medskip

  For $t_0>0$ small,
$A_{t_0}=\{x\in \Omega: \ \rho(x)<t_0\}$ is $C^2$ and let us
define
\begin{equation}\label{3.3.1}
V_\tau (x)=\left\{ \arraycolsep=1pt
\begin{array}{lll}
 \rho(x)^{ \tau},\ & x\in A_{t_0},\\[2mm]
 l(x),\ \ \ \ & x\in \Omega\setminus A_{t_0},
 \\[2mm] 0,\ &x\in\Omega^c,
\end{array}
\right.
\end{equation}
 where $\tau\in(-1,0)$ and  the function $l$  is  positive  such that $V_\tau$ is $C^2$ in $\Omega$.
\medskip

\noindent{\bf Proof of Proposition \ref{cr 1.1}.} $(i)$ Now we prove the nonexistence when $q\le 1+2\alpha$.  From Theorem 1.1 and Theorem 1.2 in \cite{CHY1},  the semilinear fracional problem
\begin{equation}\label{2.12-1}
 \arraycolsep=1pt
\begin{array}{lll}
  (-\Delta)^\alpha u+  c_1u^q=0\quad&
 {\rm in}\quad \Omega,\\[1.5mm]
  \phantom{------\ }
  u=0 & {\rm in}\quad \R^N\setminus\Omega,\\[1.5mm]
  \phantom{\ \ }
\displaystyle \lim_{x\in\Omega,x\to\partial\Omega}u(x)=+\infty
\end{array}
\end{equation}
admits a sequence solutions   $\{v_k\}_k$
 satisfying that the mapping $k\mapsto v_k$ is increasing,
$$  v_k(x)\le c_{25}k\rho(x)^{\alpha-1}, \quad \forall x\in\Omega$$
and
\begin{equation}\label{2.31}
 \lim_{k\to\infty} v_k(x)=\infty,\quad \forall x\in\Omega,
\end{equation}
 where $c_{25}>0$.

We observe that $v_k$ is a sub-solution of (\ref{eq 1.1}) for any $k$.

If (\ref{eq 1.1}) has a solution $u$ satisfying (\ref{1.3}), then by the Comparison Principle, for any $k$,
there holds that
$$
  v_k(x)\le u(x), \quad \forall x\in\Omega.
$$
Then it is impossible that $u$ is a solution of (\ref{eq 1.1}) by (\ref{2.31}).

\medskip

$(ii)$ When $q\in(1+2\alpha,\frac{1+\alpha}{1-\alpha})$, it infers from \cite{CFQ} that there exists a solution $v_q$ of (\ref{2.12-1})
replacing $c_1$ by $c_3$ from the assumption (\ref{1.5})
such that
\begin{equation}\label{2.11-1}
\frac1{c_{26}}\rho(x)^{-\frac{2\alpha}{q-1}}\le  v_q(x)\le c_{26}\rho(x)^{-\frac{2\alpha}{q-1}},\qquad\forall x\in\Omega.
\end{equation}
where $c_{26}>0$.
By (\ref{1.5}),  $v_p$ is a sub-solution of
\begin{equation}\label{eq 2.5}
 \arraycolsep=1pt
\begin{array}{lll}
  (-\Delta)^\alpha u+ f(u)=u\phi \quad&
 {\rm in}\quad \Omega,\\[1.5mm]
  \phantom{------\ }
  u=0 & {\rm in}\quad \R^N\setminus\Omega,\\[1.5mm]
  \phantom{\ \ }
\displaystyle \lim_{x\in\Omega,x\to\partial\Omega}u(x)=+\infty.
\end{array}
\end{equation}
So $v_p$ is a sub-solution of (\ref{eq 1.1}).\smallskip

We next construct a suitable super solution of (\ref{eq 1.1}). From \cite[Proposition 3.1]{CFQ}, we know that
the function $V_\tau$ with $\tau=-\frac{2\alpha}{p-1}\in(-1,0)$ satisfies
$$(-\Delta)^\alpha V_\tau(x)\ge c_\tau \rho(x)^{\tau-2\alpha},\qquad\forall x\in\Omega,$$
where $V_\tau$ is given by (\ref{3.3.1}).

We consider $\lambda V_\tau$ with $\lambda>0$. We observe that
\begin{eqnarray*}
 (-\Delta)^\alpha_\Omega (\lambda V_\tau)+ f(\lambda V_\tau) &=& (-\Delta)^\alpha  (\lambda V_\tau)+ f(\lambda V_\tau)- \lambda\phi V_\tau
 \\&\ge&
 c_\tau\lambda\rho(x)^{\tau-2\alpha}+c_2c_{26}^{-p}\lambda^p\rho(x)^{-\frac{2\alpha p}{p-1}}-c_{27}\lambda\rho(x)^{\tau-2\alpha}   \\&\ge &  \left[c_2c_{26}^{-p}\lambda^{p-1}-|c_\tau|-c_{27}\right]\lambda\rho(x)^{\tau-2\alpha}
 \\&\ge& 0
\end{eqnarray*}
if $\lambda>0$ big sufficiently.
By Theorem \ref{th 3.2.4}, it deduces that (\ref{eq 1.1}) has a solution $u$ such that
$$v_q\le u\le \lambda V_\tau\quad{\rm in}\quad \Omega,$$
which implies (\ref{1.4}).
\qquad $\Box$

\setcounter{equation}{0}
\section{ Boundary blow-up solutions for $\alpha\in(\frac12,1)$}


\subsection{Existence}
Denote by $G_{\Omega,\alpha}$ the Green kernel of $(-\Delta)^\alpha_\Omega$ in $\Omega\times\Omega$ and by $\mathbb{G}_{\Omega,\alpha}[\cdot]$ the
Green operator defined as
$$\mathbb{G}_{\Omega,\alpha}[g](x)=\int_{\Omega} G_{\Omega,\alpha}(x,y)g(y)dy. $$
\begin{proposition}\label{pr 2.2}
Assume that  $\alpha\in(\frac12,1)$, $n\in\N$ and $g\in C^{\theta}(\bar\Omega)$ with $\theta>0$,  then
\begin{equation}\label{2.4}
 \arraycolsep=1pt
\begin{array}{lll}
 \displaystyle  (-\Delta)^\alpha_\Omega w=g\qquad & {\rm in}\quad   \Omega,\\[2mm]
\phantom{ (-\Delta)^\alpha  }
 \displaystyle   w=n\quad & {\rm on}\quad   \partial  \Omega
\end{array}
\end{equation}
admits a unique solution $w_n$ such that
\begin{equation}\label{2.5}
-\mathbb{G}_{\Omega,\alpha}[g_-]\le w_n-n\le \mathbb{G}_{\Omega,\alpha}[g_+]\quad{\rm in}\quad \Omega,
\end{equation}
where $g_\pm=\max\{\pm g,0\}$.

\end{proposition}
{\bf Proof.}  {\it Existence.}  Since $\mathbb{G}_{\Omega,\alpha}[g]$ is a solution of  \arraycolsep=1pt
$$
 \displaystyle  (-\Delta)^\alpha_\Omega w=g\quad  {\rm in}\quad   \Omega,
$$
From \cite{CK}, there exists $c_{28}>0$ such that for any $(x,y)\in
\Omega\times\Omega$ with $x\neq y$,
\begin{equation}\label{annex 01}
G_{\Omega,\alpha}(x,y)\le c_{28}
\min\left\{\frac1{|x-y|^{N-2\alpha}},\frac{\rho(x)^{2\alpha-1}\rho(y)^{2\alpha-1}}{|x-y|^{N-2+2\alpha}} \right\}.
\end{equation}
 For $x\in\Omega$, we have that
\begin{eqnarray*}
 | \mathbb{G}_{\Omega,\alpha}[g](x)| &\le & c_{28}\int_\Omega \frac{\rho(x)^{2\alpha-1}\rho(y)^{2\alpha-1}}{|x-y|^{N-2+2\alpha}}  |g(y)|dy\\
    &\le & c_{28}\rho(x)^{2\alpha-1}\norm{g}_{L^\infty (\Omega)}\int_\Omega  \frac{\rho(y)^{2\alpha-1}}{|x-y|^{N-2+2\alpha}}dy
  \\& \le& c_{29}\norm{g}_{L^\infty (\Omega)}\rho(x)^{2\alpha-1},
\end{eqnarray*}
where $c_{29}>0$.
Therefore, $\mathbb{G}_{\Omega,\alpha}[g]$ is a solution of
\begin{equation}\label{eq 2.3}
\arraycolsep=1pt
\begin{array}{lll}
 \displaystyle  (-\Delta)^\alpha_\Omega w=g\qquad & {\rm in}\quad   \Omega,\\[2mm]
\phantom{ (-\Delta)^\alpha  }
 \displaystyle   w=0\quad & {\rm on}\quad   \partial  \Omega.
\end{array}
\end{equation}
and $n+\mathbb{G}_{\Omega,\alpha}[g]$ is obvious a solution of (\ref{2.4}).

{\it Uniqueness.} Let $v$ be another solution of (\ref{2.4}), we observe that $w-v$ is a solution of
$$\arraycolsep=1pt
\begin{array}{lll}
 \displaystyle  (-\Delta)^\alpha_\Omega u=0\qquad & {\rm in}\quad   \Omega,\\[2mm]
\phantom{ (-\Delta)^\alpha  }
 \displaystyle   u=0\quad & {\rm on}\quad   \partial  \Omega.
\end{array}$$
Then it follows by Maximum Principle that $w-v\equiv0$ in $\Omega$.

Finally, since $\mathbb{G}_{\Omega,\alpha}[g_+]$ is a super-solution of (\ref{eq 2.3}) and $-\mathbb{G}_{\Omega,\alpha}[g_-]$
is a sub-solution of (\ref{eq 2.3}), then (\ref{2.5}) follows.\qquad$\Box$

\smallskip

We remark that the existence of solution to (\ref{2.4}) could be extended into the one general     boundary data. Precisely,
let $\xi:\partial\Omega\to \R$ be a boundary trace of a $C^2(\bar\Omega)$ function $\tilde \xi$, i.e.
$$
\xi=\tilde \xi\quad {\rm on}\quad \partial\Omega.
$$
For  $\alpha\in(\frac12,1)$,  problem
\begin{equation}\label{re 2.4}
 \arraycolsep=1pt
\begin{array}{lll}
 \displaystyle  (-\Delta)^\alpha_\Omega w=0\qquad & {\rm in}\quad   \Omega,\\[2mm]
\phantom{ (-\Delta)^\alpha  }
 \displaystyle   w=\xi\quad & {\rm on}\quad   \partial  \Omega
\end{array}
\end{equation}
admits a unique solution
\begin{equation}\label{re 2.5}
w_\xi=\tilde \xi-\mathbb{G}_{\Omega,\alpha}[(-\Delta)^\alpha_\Omega\tilde \xi ]\quad{\rm in}\quad \Omega.
\end{equation}
 We observe that $\mathbb{G}_{\Omega,\alpha}[(-\Delta)^\alpha_\Omega\tilde \xi ]$ decays at the rate $\rho^{2\alpha-1}$ and
 $w_\xi$ is independent of the choice of $\tilde \xi$. In fact, let $\tilde \xi_1\in C^2(\bar\Omega)$ have the trace $\xi$ and
 the corresponding solution $v_\xi$
 then $w:=w_\xi-v_\xi$ is a solution of
 $$ \arraycolsep=1pt
\begin{array}{lll}
 \displaystyle  (-\Delta)^\alpha_\Omega w=0\qquad & {\rm in}\quad   \Omega,\\[2mm]
\phantom{ (-\Delta)^\alpha  }
 \displaystyle   w=0\quad & {\rm on}\quad   \partial  \Omega,
\end{array}$$
 which implies by Strong Maximum Principle that
 $$w\equiv0.$$
In the particular case that $\xi=n$, we have that $\tilde \xi=n$ in $\Omega$ and $\mathbb{G}_{\Omega,\alpha}[(-\Delta)^\alpha_\Omega\tilde \xi]=0$ in $\Omega$.

This subsection is devoted to study the existence of  solution of (\ref{eq 3.1}).
To this end, we first introduce following lemma.
\begin{lemma}\label{lm 3.2}
Let $n\in\N$, $b\ge 0$ and $g\in C^1(\bar\Omega)$,  then
\begin{equation}\label{eq 3.3}
  \arraycolsep=1pt
\begin{array}{lll}
 \displaystyle  (-\Delta)^\alpha_\Omega   u+b u=g\quad & {\rm in}\quad   \Omega,\\[2mm]
\phantom{ (-\Delta)^\alpha +bu }
 \displaystyle   u=n\quad & {\rm on}\quad   \partial  \Omega
\end{array}
\end{equation}
admits a unique solution.

\end{lemma}
{\bf Proof.}   
We observe that $n+\mathbb{G}_{\Omega,\alpha}[g_+]$ and $n-\mathbb{G}_{\Omega,\alpha}[g_-]$ are
super and sub-solutions of (\ref{eq 3.3}) respectively. We make an extension of  $n+\mathbb{G}_{\Omega,\alpha}[g_+]$ and $n-\mathbb{G}_{\Omega,\alpha}[g_-]$ by $n$ in $\R^N\setminus\Omega$
 and still denote $n+\mathbb{G}_{\Omega,\alpha}[g_+]$ and $n-\mathbb{G}_{\Omega,\alpha}[g_-]$.  Let $\Omega_t:=\{x\in\Omega: \rho(x)>t\}$ for $t\ge0$ and then
there exists $t_0>0$ such that $\Omega_t$ is $C^2$ for any $t\in[0, t_0]$, since  $\Omega$ is $C^2$.

By Perron's method, there exists a unique  solution $w_t$ of
$$
  \arraycolsep=1pt
\begin{array}{lll}
 \displaystyle  (-\Delta)^\alpha   u+(b+\phi)u=g-bn\quad & {\rm in}\quad   \Omega_t,\\[2mm]
\phantom{ (-\Delta)^\alpha +(b+\phi)u }
 \displaystyle   u=n-\mathbb{G}_{\Omega,\alpha}[g_-] \quad & {\rm in}\quad   \R^N\setminus  \Omega_t,
\end{array}
$$
where $\phi$ is defined as (\ref{con 1}). Since $t\in(0,t_0)$, $\phi$ is positive and $\phi\in C^{0,1}_{\rm loc}(\Omega_t)$,
then $w_t$ is a  solution of
$$
  \arraycolsep=1pt
\begin{array}{lll}
 \displaystyle  (-\Delta)^\alpha_\Omega   u+bu=g+bn\quad & {\rm in}\quad   \Omega_t,\\[2mm]
\phantom{ (-\Delta)^\alpha +bu}
 \displaystyle   u=n-\mathbb{G}_{\Omega,\alpha}[g_-] \quad & {\rm in}\quad   \Omega\setminus  \Omega_t
\end{array}
$$
and by Theorem \ref{comparison}, we derive that
$$n-\mathbb{G}_{\Omega,\alpha}[g_-]\le w_t\le w_{t'}\le n+\mathbb{G}_{\Omega,\alpha}[g_+]\quad{\rm for}\ 0< t'<t<t_0. $$
By Proposition \ref{pr 2.1} and Theorem \ref{stability}, the limit of $w_t$ as $t\to0$ is a classical solution of (\ref{eq 3.3}).
\qquad$\Box$

\medskip

\noindent{\bf Proof of Proposition \ref{pr 3.1}.}   {\it Existence.} Let us define
$$w_+(x)= \int_\Omega G_{\Omega,\alpha}(x,y)g_+(y)dy\quad{\rm and}\quad w_-(x)= \int_\Omega G_{\Omega,\alpha}(x,y)g_-(y)dy.$$
By (\ref{annex 01}), there exists $c_{30}>0$ such that
$$0\le w_+(x)\le c_{30}\norm{g}_{L^\infty(\Omega)}\rho(x)^{2\alpha-1},\quad x\in\Omega$$
and
$$0\le w_-(x)+f(n)\int_\Omega G_{\Omega,\alpha}(x,y)ndy\le c_{30}(\norm{g_-}_{L^\infty(\Omega)}+f(n))\rho(x)^{2\alpha-1},\quad x\in\Omega.$$
Let
 $$\bar w(x)=n-w_-(x)-f(n)\int_\Omega G_{\Omega,\alpha}(x,y)ndy$$
 and
 $$b_1=\max\{n+\norm{w_+}_{L^\infty(\Omega)}, \norm{\bar w}_{L^\infty(\Omega)}\},$$
then $\varphi(s):=(\norm{f'}_{L^\infty([-b_1,b_1])}+b_1)s-f(s)$ is increasing in $[-b_1,b_1]$.
  It follows by Lemma \ref{lm 3.2} that
  there exists a unique solution $v_m$ of
  \begin{equation}\label{eq 3.4}
  \arraycolsep=1pt
\begin{array}{lll}
 \displaystyle  (-\Delta)^\alpha_\Omega   v_m+b_2 v_m=b_2v_{m-1}-f(v_{m-1})+g\quad & {\rm in}\quad   \Omega,\\[2mm]
\phantom{ (-\Delta)^\alpha +bu_n- }
 \displaystyle   v_m=n\quad & {\rm on}\quad   \partial  \Omega,
\end{array}
\end{equation}
where $b_2=\norm{f'}_{L^\infty([-b_1,b_1])}+b_1$, $m\in\N$ and $v_0=-b_1$.
 We observe that $\{v_m\}$ is a increasing sequence and uniformly bounded in $\Omega$.
 Therefore, the limit of $\{v_m\}$ as $m\to\infty$ satisfies (\ref{eq 3.1}).

 {\it To prove (\ref{3.1}).} By direct computation, we have that
\begin{eqnarray*}
 (-\Delta)^\alpha_\Omega (n+w_+(x))+  f(n+w_+(x)) &\ge &  g_+(x)+f(n)  \ge g(x), \quad x\in\Omega
\end{eqnarray*}
and
\begin{eqnarray*}
  (-\Delta)^\alpha_\Omega \bar w(x)+  f(\bar w(x))\le  -g_-(x) -f(n)+f(n)  \le g(x),  \quad x\in\Omega
\end{eqnarray*}
thus $n+w_+$ and $n-w_--n\int_\Omega G_{\Omega,\alpha}(x,y)ndy$ are the super-solution and sub-solution of (\ref{eq 3.1}), respectively.
It infers (\ref{3.1}) by Theorem \ref{comparison}.
 \qquad$\Box$

 \smallskip

\begin{lemma}\label{lm 3.1}
Let $\tau\in(-1,0)$ and $V_{\tau}$ be defined in (\ref{3.3.1}), then
 \begin{equation}\label{3.3}
 |(-\Delta)^{\alpha}_\Omega V_\tau(x)|\leq c_{31}\rho(x)^{ \tau-2\alpha },\ \  \forall x\in \Omega,
 \end{equation}
where $c_{31}>0$.
\end{lemma}
 {\bf Proof.}
 We denote  $\tilde V_\tau=V_\tau$ in $\Omega$ and $\tilde V_\tau=0$ in $\R^N\setminus \Omega$,
 from \cite[Proposition 3.2]{CFQ},  there exists $c_{32}>1$ such that
 \begin{equation}\label{3.2}
 |(-\Delta)^{\alpha}\tilde V_\tau(x)|\leq c_{32}\rho(x)^{ \tau-2\alpha },\ \  \forall x\in \Omega.
 \end{equation}
We observe that
 $$(-\Delta)^{\alpha}_\Omega V_\tau(x)=(-\Delta)^{\alpha}\tilde V_\tau(x)- V_\tau(x)\phi(x), $$
 where $\phi$ is defined as (\ref{con 1}) and  by Lemma \ref{lm 2.2}, we have that
 \begin{eqnarray*}
\phi(x) \le  c_{12}\rho(x)^{-2\alpha},\ \  \forall x\in \Omega.
 \end{eqnarray*}
Together with (\ref{3.2}), we have that
 \begin{eqnarray*}
|(-\Delta)^{\alpha}_\Omega V_\tau(x)| &\le & |(-\Delta)^{\alpha}\tilde V_\tau(x)|+c_{12}V_\tau(x)\rho(x)^{-2\alpha} \\
     &\le &  c_{33} \rho(x)^{ \tau-2\alpha },\ \  \forall x\in \Omega.
 \end{eqnarray*}
 The proof ends.\qquad$\Box$

\smallskip

\noindent{\bf Proof of  Theorem \ref{teo 1}$(i)$.} From Proposition \ref{pr 3.1} with $g\equiv0$, there exists a unique positive solution $u_n$ of
\begin{equation}\label{eq 3.2}
  \arraycolsep=1pt
\begin{array}{lll}
 \displaystyle  (-\Delta)^\alpha_\Omega   u+h(u)=0\quad & {\rm in}\quad   \Omega,\\[2mm]
\phantom{ (-\Delta)^\alpha +f(u) }
 \displaystyle   u=n\quad & {\rm on}\quad   \partial  \Omega
\end{array}
\end{equation}
and
$$n-n^p\rho(x)^{\alpha-1}\le  u_n(x)\le n,\ \  \forall x\in \Omega.$$
By Theorem \ref{comparison}, for any $n\in\N,$
$$u_n\le u_{n+1}\quad  {\rm in}\quad   \Omega.$$

From lemma \ref{lm 3.1}, there exists $\lambda>0$ such that
$\lambda V_{-\frac{2\alpha}{p-1}}$ is a super-solution of (\ref{eq 3.2}),
where $-\frac{2\alpha}{p-1}\in(-1,0)$ for $p>1+2\alpha$.
It follows by Theorem
\ref{comparison}, we have that for all $n\in\N$,
$$u_{n}\le \lambda V_{-\frac{2\alpha}{p-1}}\quad  {\rm in}\quad   \Omega.$$
Then the limit of $\{u_n\}$ exists in $\Omega$, denoting by $u_\infty$.
Moreover, we have that $u_n$ has uniformly bound in $L^\infty$ locally in $\Omega$,
and then by regular result, we infer that  $u_n$ has uniformly bound in $C^{2\alpha+\theta}$ locally in $\Omega$.
By Theorem \ref{stability}, $u_\infty$ is a viscosity solution of (\ref{eq 1.1}).

{\it Lower bound.} From Proposition \ref{pr 3.1}, we have that
$$u_n\ge n-c_{34}n^q\rho^{2\alpha-1}\quad{\rm in}\quad \Omega,$$
then for $n$ big, let $r=(\lambda n)^{-\frac{q-1}{2\alpha-1}}$, where    $\lambda=(2^{2\alpha}c_{34})^{\frac1{q-1}}$ chosen later,
then for $x\in \Omega_r\setminus \Omega_{2r}$,  we have that
\begin{eqnarray*}
 u_n(x) &\ge & \frac1\lambda r^{-\frac{2\alpha-1}{q-1}} -c_{34}\frac1{\lambda^p} r^{-\frac{2\alpha-1}{q-1}p}(2r)^{2\alpha-1} \\
    &\ge&\frac1\lambda(1-\frac{2^{2\alpha-1}c_{34}}{\lambda^{q-1}}) r^{-\frac{2\alpha-1}{q-1}}\\
    &\ge & \frac1{2\lambda} \rho(x)^{-\frac{2\alpha-1}{q-1}}.
\end{eqnarray*}
where $\lambda$ is independent of $n$. For any $x\in \Omega\setminus \Omega_{r_0}$, there exists $n$ such that
$$u_\infty(x)\ge u_n(x)\ge \frac1{2\lambda} \rho(x)^{-\frac{2\alpha-1}{q-1}}. $$

We notice that the solution $u_\infty$ is the minimal solution of (\ref{eq 1.1}), since for any boundary blow-up solution $u$,
we may imply by   Comparison Principle that $u\ge u_n$  in $\Omega,$
which infers that $u_\infty\le u$ in $\Omega$.
The proof ends. \qquad $\Box$
\smallskip

\subsection{Nonexistence }

This subsection is devoted to prove the nonexistence part of Theorem \ref{teo 1}. \smallskip

\noindent{\bf Proof of Theorem \ref{teo 1} $(ii)$.} If $q\le1$,
we observe that for $n>1$,
$$u_n\ge nu_1\quad {\rm in}\quad \Omega,$$
which implies that (\ref{eq 1.1}) has no solution.

In what follows, we assume that $q>1$.
 By contradiction, we may assume that there exists a solution $u$ of (\ref{eq 1.1})
when $f(s)\le c_{11}s^q$  for $s\ge0$ and  $q$ satisfying (\ref{q+1}).
 By Theorem \ref{comparison}, we have that
$$u_n\le u\quad{\rm in}\quad \Omega.$$
From Proposition \ref{pr 3.1}, we have that
$$
 u_n\ge n-c_{34}n^q\rho^{2\alpha-1}\quad{\rm in}\quad \Omega.
$$
Then for $n$ big, let $r_n=(\lambda n)^{-\frac{q-1}{2\alpha-1}}$, where    $\lambda=(2^{2\alpha}c_{34})^{\frac1{q-1}}$ chosen later,
then for $x\in \Omega_{r_n}\setminus \Omega_{2r_n}$,  we have that
\begin{eqnarray*}
 u_n(x)  \ge  \frac1\lambda r_n^{-\frac{2\alpha-1}{q-1}} -\frac{c_{34}}{\lambda^p} r_n^{-\frac{2\alpha-1}{q-1}p}(2r_n)^{2\alpha-1}
     \ge   \frac1{2\lambda} \rho(x)^{-\frac{2\alpha-1}{q-1}}.
\end{eqnarray*}
 For any $x\in \Omega\setminus \Omega_{r_0}$, there exists $n$ such that
\begin{equation}\label{4.1}
 u(x)\ge u_n(x)\ge \frac1{2\lambda} \rho(x)^{-\frac{2\alpha-1}{q-1}}.
\end{equation}

When $1<q\le 2\alpha$, we have that $\rho^{-\frac{2\alpha-1}{q-1}}$ is not in $L^1(\Omega)$,
then it follows from (\ref{4.1}) for any $x\in\Omega$ and any $\epsilon>0$
\begin{eqnarray*}
 (-\Delta)^\alpha_{\Omega,\epsilon} u(x)&\le & -\int_{\Omega\setminus B_\epsilon(0)} \frac{u_n(y)-u(x)}{|x-y|^{N+2\alpha}}dy \\
    &\le& -\epsilon^{-N-2\alpha}\left[\int_{\Omega} u_n(y) dy- u(x)|\Omega|\right]
    \\&\to&-\infty\quad{\rm as}\quad n\to\infty,
\end{eqnarray*}
which is impossible.
\smallskip

From (\ref{q+1}), we have that $-\frac{2\alpha-1}{q-1}<\alpha-1$, then if follows from (\ref{4.1}) that
\begin{equation}\label{4.2}
\lim_{\rho(x)\to0^+}u(x)\rho^{1-\alpha}(x)=+\infty,
\end{equation}
which contradicts Proposition \ref{cr 1.1} $(i)$. \qquad$\Box$

\medskip

\noindent{\bf Acknowledgements:} H. Chen is supported by NSFC of China grant 11401270
and
the Project-sponsored by SRF for ROCS, SEM.

Huyuan Chen
\medskip

\noindent Department of Mathematics, Jiangxi Normal University,

\noindent Nanchang, Jiangxi 330022, PR China

  and

\noindent  Institute of Mathematical Sciences,  New York University  Shanghai,

\noindent Shanghai 200120, PR China

 \bigskip

 Hichem Hajaiej

 \medskip
\noindent  Institute of Mathematical Sciences,  New York University  Shanghai,

 \noindent Shanghai 200120, PR China


\begin{thebibliography}{99}
\bibitem {AGQ} S. Alarc\'{o}n, J. Garc\'{i}a-Meli\'{a}n and A. Quaas, Keller-Osserman type conditions for some elliptic
problems with gradient terms, {\it J. Diff. Eq. 252}, 886-914 (2012).

\bibitem {BM1}
C. Bandle and M. Marcus, Large solutions of semilinear elliptic
equations: Existence, uniqueness and asymptotic behaviour, {\it J.
Anal. Math. 58,} 9-24 (1992).

\bibitem {BM} C. Bandle and M. Marcus,
Asymptotic behaviour of solutions and derivatives for semilinear elliptic problems with blow-up on the boundary,
{\it  Ann. Inst. H. Poincar\'{e} Anal. Non Lin\'{e}aire 12,} 155-171 (1995).

\bibitem {BBC}
K. Bogdan, K. Burdzy and Z. Chen,   Censored stable processes,
 {\it Probability theory and related fields 127,} 89-152 (2003).

\bibitem {CS}  L. Caffarelli   and L. Silvestre,
Regularity theory for fully nonlinear integro-differential equations,
{\it  Comm.  Pure  Appl. Math. 62},  597-638 (2009).


\bibitem {CFQ} H. Chen, P. Felmer and A. Quaas, Large solution to elliptic  equations involving fractional Laplacian,
 {\it  Ann. Inst. H. Poincar\'{e} Anal. Non Lin\'{e}aire}, to appear, DOI: 10.1016/j.anihpc.2014.08.001.

\bibitem {CV3} H. Chen and L. V\'{e}ron,
Weakly and strongly singular solutions of
semilinear fractional elliptic equations, {\it Asym. Anal. 88,} 165-184 (2014).

\bibitem {CHY1} H. Chen, H. Hajaiej and Y. Wang,  Boundary blow-up solutions to fractional elliptic
equations in a measure framework, {\it  arXiv:1505.02490}  (2015).

\bibitem {CK} Z. Chen and P. Kim, Green function estimate for censored stable processes,
{\it Probability theory and related fields 124,} 595-610 (2002).

\bibitem {CKS} Z. Chen, P. Kim and R.Song,  Two-sided heat kernel estimates for censored stable-like processes,
 {\it Probability theory and related fields 146,} 361-399 (2010).

\bibitem{Duan2015} J. Duan, An Introduction to Stochastic Dynamics, {\it Cambridge University Press, New York} (2015).


\bibitem {DF}
B. Dyda and  R. Frank, Fractional Hardy-Sobolev-Maz'ya inequality for domains,  {\it Studia Math. 208, }  151-166  (2012).

\bibitem {DZZ} Y. Du, Z. Guo and F. Zhou, Boundary blow-up solutions with interior layers and spikes in a bistable problem,
{\it Discrete Contin. Dyn. Syst. 19,} 271-298 (2007).

\bibitem {FQ0} P.  Felmer and A. Quaas, Boundary blow-up solutions for fractional elliptic equations, {\it
 Asym. Anal. 78,}  123-144 (2012).


\bibitem {G}
Q. Guan,   Integration by parts formula for regional fractional laplacian,
{\it Comm. Math. Phys.  266,} 289-329 (2006).

\bibitem {GM}
Q. Guan and  Z. Ma,  Reflected symmetric ¦Á-stable processes and regional fractional Laplacian,
{\it Probability theory and related fields 134,} 649-694 (2006).


\bibitem {K} J. B. Keller, On solutions of $\Delta u = f(u)$, {\it Comm. Pure Appl.
Math.  10,}  503-510 (1957).

\bibitem {LP} T. Leonori and A. Porretta, The boundary behavior of blow-up solutions related to a stochastic control problem with state constraint,
{\it SIAM J. Math. Anal. 39,} 1295-1327 (2008).


\bibitem {MV0} M. Marcus  and L. V\'{e}ron, Uniqueness and asymptotic behavior of solutions with boundary blow-up for a class of nonlinear elliptic equations,
{\it  Ann. Inst. H. Poincar\'{e} Anal.  Non Lin\'{e}aire 14,}   237-274 (1997).



\bibitem {MV} M. Marcus and L. V\'{e}ron, Existence and uniqueness results for large solutions of
general nonlinear elliptic equation, {\it J. Evol. Eq. 3,} 637-652 (2003).

\bibitem {O} R. Osserman, On the inequality $\Delta u = f(u)$, {\it Pac. J. Math. 7},
1641-1647 (1957).

\bibitem {TR}
E. Topp and J. Rossi,  Large Solutions for a Class of Semilinear Integro-Differential Equations with Censored Jumps,
{\it ACTA DE RESUMENES LXXXIII Encuentro Anual Sociedad de Matem¨¢tica de Chile, 115.}

\bibitem {RS} X. Ros-Oton and J. Serra, The Dirichlet problem for the fractional
laplacian: regularity up to the boundary, {\it J. Math. Pures Appl. 101},  275-302 (2014).



\end{thebibliography}
\end{document}